\documentclass[11pt,reqno]{amsart}
\usepackage{amssymb,latexsym}
\usepackage{graphicx}
\usepackage{fancyhdr,amssymb}
\usepackage{url}		

\theoremstyle{plain}
\numberwithin{equation}{section}

\begin{document}
\fancyhead{}
\renewcommand{\headrulewidth}{0pt}
\fancyfoot{}
\fancyfoot[LE,RO]{\medskip \thepage}

\setcounter{page}{1}

\title[Generalised Binomial coefficients and Jarden's Theorem]{Generalised Binomial coefficients
 and Jarden's Theorem}
\author{Cheng Lien Lang}
\address{Department Applied of Mathematics\\
                I-Shou University\\
                Kaohsiung, Taiwan\\
                Republic of China}
\email{cllang@isu.edu.tw}
\thanks{}
\author{Mong Lung Lang}
\address{  Singapore 669608, Republic of Singapore   }
\email{lang2to46@gmail.com}

\begin{abstract}
We prove a stronger version of  Jarden's Theorem on recurrence for product of recursive functions.

\end{abstract}

\maketitle

\vspace{-.8cm}

\section{Generalised Binomial Coefficients}

\noindent Let $\{x_0, x_1, x_2, \cdots\}$ be a sequence defined by the following three term recurrence
($x_0 = a, x_1= b$ for some $a, b,p,q \in \Bbb C$)
 $$x_r= px_{r-1} -qx_{r-2}. \eqno(1.1)$$

 \noindent
 In the case $x_0=0$ and $x_1=1$, we shall denote this sequence by $\{u_r\}$. Equivalently, $\{u_r\}$ is given by
  $$u_0=0, u_1=1, u_r= pu_{r-1}-qu_{r-2}.\eqno(1.2)$$

 \medskip
 \noindent
  Let $\sigma$ and $ \tau$ be roots of $x^2 -px+q=0$ (known as the characteristic polynomial
  of (1.1)). It is well known that
 $$u_r = \sum_{i=1}^{r-1} \sigma^{r-1-i}\tau^i = f_r(\sigma, \tau) ,\eqno(1.3)$$

 \medskip
\noindent where $f_r(x,y)$ is the polynomial  $ \sum x^{r-1-i}y^i  =(x^r-y^r)/(x-y).$
 Let $G(z)$ be the following.
 $$ G(z) = (1-z^{m+n})(1-z^{m+n-1})\cdots (1-z^{m+1}/(1-z^n)(1-z^{n-1})\cdots (1-z).\eqno(1.4)$$
 $G(z)$ is known as the Gaussian binomial coefficient  and is a polynomial in $z$ (a more powerful result actually implies that (1.4) can be written
  as product of cyclotomic polynomials).
  One may apply this fact to show that
  that
 $(f_r f_{r-1}\cdots f_{r-k+1})/(f_kf_{k-1}\cdots f_1)\in \Bbb Z[x,y]$ is  a polynomial in $x$ and $y$.
Denoted by $F(r,k, x, y)$ this polynomial. Define

 $$ (r|k)_u = F(r,k,\sigma, \tau),\eqno(1.5)$$

\medskip
\noindent for $1\le k\le r$ . In the case $k=0$, we define $(r|0)_u =1$.
Note that $p$ and $q$ are arbitrary.
It is clear that if $u_1u_2\cdots u_k\ne0$, then

$$(r|k)_u =   u_ru_{r-1}\cdots u_{r-k+1}/u_ku_{k-1}\cdots u_1,\eqno(1.6)$$

\medskip
\noindent
where the right hand side of (1.6) is the
 generalised binomial coefficients of $\{u_r\}$.
 Note that $(r|k)_u$ is well defined for any sequence $\{u_r\}$ but  $u_ru_{r-1}\cdots u_{r-k+1}/u_ku_{k-1}\cdots u_1$
  is well defined only if  $u_ku_{k-1}\cdots u_1\ne0$.
Hence, we may treat $(r|k)_u $ as a generalisation of the generalised
 binomial coefficients of $\{u_r\}$. However, as identity (1.6) suggested, we propose to call
  the  more general expression
  $(r|k)_u$ the {\em generalised binomial coefficient} as well.
   To ensure this  is not just another generalisation of no significance, we
     will prove, in the following section,  a stronger version of Jarden's Theorem [6] which  plays a very important role in the
      study of recurrence relations.
         We shall also take this opportunity to report (in  Appendix $A$)  our generalisations
          as well as an alternative proof (which uses $(r|k)_u$) of a result
    by Bachmann [1], Carmichael [3], and Jarden and Motzkin (see [6]).
    Application of Jarden's Theorem can be found in Section 3.

\section {Recurrence for product  of recurrence functions}

\medskip
\noindent {\bf Theorem  2.1.} (Jarden) {\em  Let $n \in \Bbb N$ be fixed and let $X(m)$ be a product
 of $n$ functions, each of which satisfies $(1.1)$.
  Suppose that $u_1u_2\cdots u_{n+1}\ne 0$.
 Then $X(m)$ satisfies the
  following recurrence.
   $$
  \sum_{i=0}^{n+1}
  (-1)^{i} q^{i(i-1)/2} \frac{u_{n+1}u_n\cdots u_{n-i+2}}{u_i u_{i-1}\cdots u_1}
 X(m-i) =0.\eqno(2.1)$$}

 \noindent {\em The polynomial  $
  \sum_{i=0}^{n+1}
  (-1)^{i} q^{i(i-1)/2}
  (u_{n+1}u_n\cdots u_{n-i+2})x^i/(u_i u_{i-1}\cdots u_1)$
   is called the  characteristic polynomial associated to $X(m)$.
 }

\medskip
\noindent {\em Proof.} See [4], [6]. \qed

\medskip
 With our generalisation (1.5), the assumption
 $u_1u_2\cdots u_{n+1}\ne 0$ in Theorem 2.1 can be removed. As a consequence, we  have the following
  stronger version of Jarden's Theorem. Note that Theorem 2.2 works for all sequences $\{u_m\}$ while
   Theorem 2.1 fails for sequence such as $u_n = pu_{n-1}-qu_{n-2}$, where
   $p^4-4p^2q+3q^2=0$ $(u_6=0$).

 \medskip
 To the best of our knowledge, the known proof of Theorem 2.1 in the literature requires the fact
 that $u_1u_2\cdots u_{n+1}\ne 0$
  Consequently, the proof of Theorem 2.1 (in the literature) cannot
  be used directly without adjustment to prove Theorem 2.2.

\medskip
\noindent {\bf Theorem 2.2.} (Jarden)  {\em  Let $n \in \Bbb N$ be fixed and let $X(m)$ be a product
 of $n$ functions, each of which satisfies $(1.1)$.
 Then $X(m)$ satisfies the
  following recurrence.
   $$
  \sum_{i=0}^{n+1}
  (-1)^{i} q^{i(i-1)/2}  {(n+1|i)_u}
 X(m-i) =0,\eqno(2.2)$$}

\noindent {\em Proof.} Similarly to (1.1), we define the sequence $\{v_r\}$ by
 $v_0=0, v_1=1$,
 $$v_r= pv_{r-1}-zv_{r-2},\eqno(2.3)$$

 \medskip
 \noindent where $z$ is a variable. One sees easily that $v_r$ is a polynomial in $z$ for all $r$. We now define
  the following rational function (in $z$)

  $$  C(z, i) =\frac{v_{n+1}v_n\cdots v_{n-i+2}}{v_i v_{i-1}\cdots v_1}.\eqno(2.4)$$

 \medskip
 \noindent
 Since $v_iv_{i-2}\cdots v_1$ is a polynomial in $z$, there exists an $\epsilon_0$ such that
  $v_iv_{i-1}\cdots v_1 \ne 0$ as long as $z \in I_q = (q-\epsilon_0,q+\epsilon_0)-\{q\}.$
  Let $\{s_j\} \subset I_q$ be a sequence converges to $q$
 and let $X_j(m)$ be a product of $n$ functions, each of which satisfies (2.3) (with $z = s_j$).
 Further, the initial conditions for these $n$ functions (whose product is $X_j(m)$) are the same as the initial conditions for those $n$ functions (whose product is $X(m)$).
  Applying Theorem 2.1, the characteristic polynomial of $X_j(m)$ is given by

   $$\Phi_{ps_j}(x) =
  \sum_{i=0}^{n+1}
  (-1)^{i} s_j^{i(i-1)/2}C(s_j, i)
 x^i.\eqno(2.5)$$

 \medskip
 \noindent
 Let $m$ be given.
 Since $\{s_j\} \to q$,
 one must have
   $\{C(s_j, i)\} \to
   (n+1|i)_u$ (see $(A2)$ of  Appendix $A$) and
  $\{X_j(w)\} \to X(w)$. Hence for every  $\epsilon \in (0, 1/2)$, there exists some $t$
   such that if $j\ge t$, then

   $$ |s_j^{i(i-1)/2} - q^{i(i-1)/2}| < \epsilon,\,\,|X(m-i) - X_j(m-i)| < \epsilon,  \,\,\,|C(s_j,i)-
   (n+1|i)_u| < \epsilon ,\eqno (2.6)$$

 \medskip
\noindent
for every $i$, where $1\le i \le n+1$ (note that $n$ is fixed).
It is clear that $\{X(m-i) \,:\, 1\le i \le n+1\}$ is bounded.
 By (1.5),
 $(n+1|i)_i$ is a polynomial in $\sigma$ and $\tau$. Hence
 $\{(n+1|i)_u \,:\, 1\le i \le n+1\}$    is  bounded. Hence
 there exists $M > 1$ such that for all $j \ge t$, $1\le i\le n+1$,

 $$|\,s_j^{i(i-i)/2} \,| < M,\,\,\,| X_j(m-i)| <M,\,\,\,|C(s_j, i)|<M.\eqno(2.7)$$
\noindent   Let
  $$ A=
 \sum_{i=0}^{n+1}
  (-1)^{i} q^{i(i-1)/2}(n+1|i)_uX(m-i).\eqno(2.8)$$

 \medskip
 \noindent
  Applying the three inequalities we have in (2.6), we may rewrite $A$ in  (2.8) into the following.

{\small $$
 A =
  \sum_{i=0}^{n+1}(-1)^i(s_j^{i(i-1)/2} +\epsilon_{1i})
  (X_j(m-i)+ \epsilon_{2i} ) (C(s_j, i) +\epsilon_{3i}).\eqno(2.9)$$}

 \medskip
 \noindent where $0<|\epsilon _{ki}| < \epsilon <1.$  Since $\Phi_{ps_j}(x)$ is the characteristic polynomial
  of $X_j(m)$ (see (2.5)), one has
 $$
  \sum_{i=0}^{n+1}
  (-1)^{i} {s_j}^{i(i-1)/2}C(s_j, i)X_j(m-i) = 0. \eqno(2.10)$$

  \medskip
  \noindent Multiply through the terms to the right  of (2.9) (which gives all together eight terms) and
  applying (2.7), (2.9), (2.10) and the fact that $0<|\epsilon _{ki} |< \epsilon <1$, $M>1$, we conclude that

$$ |A|= \big |\,\sum
  (-1)^{i} q^{i(i-1)/2}(n+1|i)_uX(m-i)\,\big | \le   0+ 6(n+2)M^2 \epsilon
  +(n+2)\epsilon
  .\eqno(2.11)$$

  \medskip
  \noindent  Hence $A=0$.  This implies that $X(m)$ admits a recurrence relation
   described by the polynomial $\Phi_{pq}(x) =\sum_{i=0}^{n+1}
  (-1)^{i} q^{i(i-1)/2}(n+1|i)_ux^i.$\qed

  \medskip

\noindent {\bf Proposition 2.3.} {\em  Let $n \in \Bbb N$ be fixed and let $X(m)$ be a product
 of $n$ functions, each of which satisfies $(1.1)$.
 Suppose that $u_k= 0$, for some $k$, where $1\le k\le n+1$.
 Then $X(m)$ satisfies the
  following recurrence.
   $$X(k+r+1)= u_{k+1}^nX(r+1)
   .\eqno(2.12)$$
  }

\noindent {\em Proof.} Suppose that $x_m$ satisfies (1.1). One can show easily by mathematical induction that
(Theorem 2.1 of [5])

$$x_{m+r+1}= u_{m+1}x_{r+1}-qu_mx_r.\eqno(2.13)$$

\medskip
\noindent Since $u_{k}=0$,
Identity (2.13) now becomes $x_{k+r+1} = u_{k+1}x_{r+1}.$
 Since $X(m)$ is a product of $n$ functions, each of which satisfies  $x_{k+r+1} = u_{k+1}x_{r+1}$,
  one has
  $X(k+r+1) = u_{k+1}^nX(r+1)$.
  Since this holds  for every  $r \in \Bbb Z$ and $k$ is fixed, one has just obtained
   a recurrence for $X(m).$\qed

\medskip

  \noindent {\bf Example 2.4.} Let $u_0=0, u_1 =1$ and $u_n= 2u_{n-1} -4u_{n-2}$. Then
  $u_2 = 2, u_3 =0, u_4 = -8$. By Proposition 2.3 and Theorem 2.2,
   $ u_m^3$
   satisfies the following relations
 $$ u_{r+4}^3 + 8^3u_{r+1}^3=0,\,\,
u_{r+4}^3 -8u_{r+3}^3 -512u_{r+1}^3 +4096u_r^3=0.\eqno (2.14)$$

\noindent Let $u_0=0, u_1 =1$ and $u_n= u_{n-1} -u_{n-2}$. Then $(7|3)_u=2$ and $X(m)=u_m^6$
 satisfies the relation
 $$X(m)-X(m-1)-2X(m-3)+2X(m-4) +X(m-6)-X(m-7)=0.\eqno(2.15)$$

\medskip
\noindent {\bf Discussion 2.5.}
 A setback for Theorem 2.2 is that $(r|k)_u$ is not
    easy to determine by its definition (1.5). We will show in Appendix A that this can be saved
     (see $(A2)$).

\section{application}

Jarden's Theorem makes the verification of many identities {\em trivial}. Take H-635 of [10]
for instance, the readers are asked to prove the following identity.
$$
\left |
\begin{array}{ccc}
F_{n} & F_{n+1}& F_{n+2}\\
F_{n+2} & F_{n}& F_{n+1}\\
F_{n+1} & F_{n+2}& F_{n}\\
\end{array}
\right | = 2(F_{n}^3 +F_{n+1}^3).\eqno(3.1)$$

The identity is proved by J. A. Sellers [10]. Alternatively,
one may apply Jarden's Theorem and conclude
 that both the left and right hand side of (3.1) satisfy the same recurrence.
 $$ X(n) = 3X(n-1) +6X(n-2)-3X(n-3) -X(n-4).\eqno(3.2)$$
  Denoted by $A(n)$ and $B(n)$ the left and right hand side of (3.1) respectively.
   Since $A(n)$ and $B(n)$ satisfy the same recurrence, one has $A(n) = B(n)$ if and
    only if $A(n) = B(n)$ for $n = 1,2,3,4$ which makes the verification of (3.1)
     straightforward.
     The same technique can be applied to
     B-1016, 1022, 1029, 1032, 1044, 1049, 1058, 1059, 1071, 1077, 1080, 1081, 1112, 1118, 1127,
         H-619, H-631, H-640 of the journal {\em Fibonacci quarterly} and many others (see [7], [8]).
        Note that the following elementary facts is needed while various identities are verified.

   \begin{enumerate}
   \item[(i)] Let $k, r\in \Bbb Z$ be fixed. Then $A(n)$ and $ kA(n+r)$ satisfy the
    same recurrence.
    \item[(ii)] Suppose that $A(n)$ and $B(n)$ satisfy the same recurrence. Then
    $A(n) \pm B(n)$ satisfies the same recurrence.
    \end{enumerate}

Let $k\in \Bbb N$ be fixed and
let $x_n$ and $u_n$ be given as in (1.1) and (1.2).
  Applying (2.13),
(i) and (ii) of the above and induction,  one has

     \begin{enumerate}
   \item[(iii)] $x_{kn}$ and $x_n^k$ satisfy the same recurrence.
   In particular, $ A(n) = F_{kn}$, $B(n) = F_n^k$, $C (n) = L_{kn}$ and $D(n) = L_n^k$
    satisfy the same recurrence.
    \end{enumerate}

\section { Appendix A : Properties of $(n|k)_u$}
\noindent Recall that $v_r$ is defined by $v_0=0, v_1=1$ and $v_r=pv_{r-1}-zv_{r-2}$.
Let $\sigma_z$ and $\tau_z$ be roots of $x^2 - px +zx=0$. Applying (1.3), (2.3) and (2.4)

$$C(z, i)  =\frac{v_{n+1}v_n\cdots v_{n-i+2}}{v_i v_{i-1}\cdots v_1} =
\frac{f_{n+1}f_n\cdots f_{n-i+2}}{f_i f_{i-1}\cdots f_1}
=F(n+1, i, x,y)\,\big |\,_{x=\sigma_z, y = \tau_z}.
\eqno(A1)
$$

\medskip
\noindent {\bf Lemma A1.} {\em $C(z, i) \in\Bbb C[z]$ is a polynomial in $z$.}

\medskip
\noindent {\em Proof.}
 Since $C(z,i) = v_{n+1}v_n\cdots v_{n-i+2}/v_i v_{i-1}\cdots v_k$,
 $C(z,i)$ is a rational function in $z$. By $(A1)$,
 $C(z,i) =F(n+1, i, \sigma_z, \tau_z)$ is
 a polynomial in $\sigma_z, \tau_z$.  Hence $C(z,i)$  must be a polynomial in $z$.\qed

\medskip
\noindent
 It is clear that if $z\to q$, then $\sigma_z\to \sigma$ and $\tau_z\to \tau$. Applying (1.5) and $(A1)$ of the
  above, one has
$$\lim_{z\to q} C(z, i) = (n+1|i)_u.\eqno (A2)$$

\medskip
\noindent
 Bachmann [1], Carmichael [3], and  Jarden and Motzkin [6] have shown that
  the generalised binomial coefficient
 $x_nx_{n-1}\cdots x_{n+1-k}/x_kx_{k-1}\cdots x_1$
  is an integer, where
  $x_0 =0, x_1 =1$ and $ x_m= ax_{m-1} +bx_{m-2},  gcd\,(a,b)=1$. Lemma $A.1$
   can be viewed as an analogue of their result. The following  is a
     generalisation of their result.

    \medskip
    \noindent {\bf Proposition A2.} {\em Let $u_r$ be given as in $(1.2)$. Suppose that
     $p, q \in \Bbb Z$. Then $(n+1|i)_u\in \Bbb Z$. }

\medskip
\noindent {\em Proof.}
By(1.5), $(n+1|i)_u$ is an algebraic integer. Since $p, q\in \Bbb Z$,
 by Lemma $A.1$, $C(z,i) \in \Bbb Q[z]$. Hence
$\lim_{z\to q} C(z, i) = (n+1|i)_u$ is a rational number. It follows that
 $(n+1|i)_u$ is an integer.\qed


\medskip
\section{Appendix B : Characteristic Polynomial }

\noindent Let (see Section 3 of [4]) $A_{n+1}$ be an $(n+1)\times (n+1)$ matrix whose $(r+1, c+1)$ entry is given by
 $$
\left(\begin{array}{c}
r\\
n-c\\
\end{array}\right ) p^{r+c-n}(-q)^{n-c}
\,\,\,.\eqno(B1)$$

\medskip
\noindent
 One sees easily that $A_n = E Q_n^tE^{-1}$, where $E$ is the matrix with ones on the
 counter diagonal and zeros elsewhere and $Q_n^t$ is the transpose of the following matrix
(to save space, we denote the binomial coefficient
{\tiny $
  \left(\begin{array}{c}
n\\
k\\
\end{array}\right )$} by $(n|k)$).

\medskip
{\small $$\,\,\,\, Q_n=
 \left [
 \begin{array}{llcrr}
 (n-1|0)p^{n-1}&    (n-2|0)    p^{n-2}& \cdots     &  p  &     1            \\
 (n-1|1) p^{n-2} (-q)& (n-2|1)    p^{n-3}(-q)& \cdots     & -q  &  0    \\
 (n-1|2) p^{n-3} (-q)^2   &(n-2|2)p^{n-4}(-q)^2& \cdots     & 0    & 0        \\
  \hspace{1cm}\vdots                  &      \hspace{1cm}     \vdots            & \ddots         & \vdots     & \vdots\\
  (n-1|n-2)   p (-q)^{n-2}& (n-2|n-2)(-q)^{n-2}&\cdots &0  & 0               \\
  (n-1|n-1) (-q)^{n-1}& 0&\cdots &0 & 0            \\
   \end{array} \right ].
    \eqno(B2))
   $$}

\medskip
\noindent These two matrices $A_{n+1}$ and $Q_n$ are quite useful in the study  of the recurrence of  $
 X(m)$ (see for examples, (2.6) of  [2], [4], [9]). Since $A_n$ and $Q_n$ are  similar to each other,
  they admit the same
  characteristic polynomial. The characteristic polynomial (see for examples,
  [2], [4], [9]), when $u_1 u_2 \cdots u_{n}\ne 0$,
   is given by

  $$
  \sum_{i=0}^{n}
  (-1)^{i} q^{i(i-1)/2} \frac{\,u_nu_{n-1} \cdots u_{n-i+1}\,}
  {u_1u_2\cdots u_i} x^i. \eqno(B3)$$

\medskip
\noindent This left the characteristic polynomial undecided when $u_i= 0$ for some $i$, and once again,
 the function we defined in (1.5) can be used to determine the characteristic polynomial.
 Following our technique given in the proof of Theorem 2.2, one can show that for
  any $\{u_r\}$, the characteristic
  polynomial of $A_n$ (hence $Q_n)$ is given by

 $$
  \sum_{i=0}^{n}
  (-1)^{i} q^{i(i-1)/2}(n|i)_u x^i. \eqno(B4)$$

\bigskip

\noindent MSC2010 : 11B39, 11B83.

\medskip

\end{document}